\documentclass[12pt]{article}
\usepackage{amssymb}
\usepackage{amsfonts}
\usepackage{mathrsfs}
\usepackage{graphicx}
\usepackage{amsmath}
\usepackage{epsfig}
\usepackage{verbatim}
\usepackage{url}
\usepackage{color}

\usepackage[top=1in, bottom=1in, left=.75in, right=.75in]{geometry}

\newcommand{\lab}[1]{\label{#1}}                % hides labels
\newcommand{\remove}[1]{}

\newcommand{\be}{\begin{equation}}
\newcommand{\ee}{\end{equation}}
\newcommand{\bea}{\begin{eqnarray}}
\newcommand{\eea}{\end{eqnarray}}
\newcommand{\bean}{\begin{eqnarray*}}
\newcommand{\eean}{\end{eqnarray*}}

\newtheorem{thm}{Theorem}%[section]

\newtheorem{lemma}[thm]{Lemma}

\newtheorem{prop}[thm]{Proposition}

\def\proof{\noindent{\bf Proof.\ }  }
\def\qed{~~\vrule height8pt width4pt depth0pt}

\def\ex{{\bf E}}
\def\eps{\epsilon}

\def\H{{\mathcal H}}
\def\no{\noindent}

\catcode`@=11 \@addtoreset{equation}{section}

\catcode`@=12
\date{}

\title{Analysis of the parallel peeling algorithm: a short proof}
\author{Pu Gao\footnote{Research supported by an NSERC Postdoctoral Fellowship.}\\ University of Toronto\\
pu.gao@utoronto.ca}
\begin{document}
\maketitle

Given a (hyper)graph $H$ and a positive integer $k$, the parallel peeling algorithm repeatedly removes all vertices of degree less than $k$ and their incident edges. When the algorithm terminates, the output is the $k$-core of $H$. Let $s(H)$ denote the number of rounds the algorithm takes. It was first proved by Achlioptas and Molloy~\cite{amxor} that, if $\H_r(n,p)$ is a random $r$-uniform hypergraph on $[n]$ with edge density $p=c/n^{r-1}$, where $c>0$ is a constant not equal to $c_{r,k}$, the emergence threshold of a non-empty $k$-core, then $s(H)=O(\log n)$ (here $r,k$ are both at least 2 and are not both equal to 2). Recently, a paper by Jiang, Mitzenmacher and Thaler~\cite{jmt} improved this result by showing that, if $c>c_{r,k}$, then $s(\H_r(n,c/n^{r-1}))=\Omega(\log n)$, i.e.\ the upper bound in~\cite{amxor} is tight; if $c<c_{r,k}$, then $s(\H_r(n,c/n^{r-1}))\le a_{r,k}\log\log n+O(1)$ where $a_{r,k}=1/\log((r-1)(k-1))$, which significantly improves~\cite{amxor}. The lower bound in the supercritical case is relatively easier whereas most of the technical proof of~\cite{jmt} was for the upper bound in the subcritical case. In this note, I give a very short proof of asymptotically the same upper bound as in~\cite{jmt} (with a slightly larger coefficient than $a_{r,k}$) in the subcritical case. In fact, my proof mainly combines several well-known results in literature.
I will prove the following.
\begin{thm}\lab{thm:main}
Assume $k,r\ge 2$, $(k,r)\neq (2,2)$ and $c<c_{r,k}$. Then a.a.s.\ $s(\H_r(n,c/n^{r-1}))\le (a^*_{r,k}+o(1))\log\log n$, where $a^*_{r,k}=1/\log(k(r-1)/r)$. 
\end{thm}

Here is the key lemma I use.

\begin{lemma}\lab{lem:small} Assume $k,r\ge 2$, $(k,r)\neq (2,2)$ and $c=O(1)$.
A.a.s.\ every subgraph of $\H_r(n,c/n^{r-1})$ with less than $\log^2 n$ vertices has average degree less than $r/(r-1)+\eps$ for every constant $\eps>0$.
\end{lemma}

\proof
Let $X_{s,t}$ denote the number of subgraphs of $\H_r(n,c/n^{r-1})$ with $s$ vertices and at least $t$ edges.
Then,
\[
\ex X_{s,t}\le \binom{n}{s}\binom{s^r}{t}\left(\frac{c}{n^{r-1}}\right)^t.
\]
Fix a constant $0<\eps<1$;
let $t=(1+\eps)s/(r-1)$; then
\[
\ex X_{s,t}\le \left(\frac{en}{s}\left(\frac{e(r-1)s^r}{(1+\eps)s}\frac{c}{n^{r-1}}\right)^{(1+\eps)/(r-1)}\right)^s\le\left(C\left(\frac{s}{n}\right)^{\eps}\right)^{s},
\]
for some constant $C>0$ depending only on $r$, $k$ and $c$.
Now immediately we have $\sum_{1\le s\le \log^2 n} \ex(X_{s,t})=o(1)$ and the lemma follows as each edge contributes $r$ to the total degree of a subgraph and $\eps>0$ is arbitrary. \qed

The following proposition is from~\cite[Section 8]{amxor}.

\begin{prop}\lab{prop}
Assume $k,r\ge 2$, $(k,r)\neq (2,2)$ and $c<c_{r,k}$; let $H=\H_r(n,c/n^{r-1})$. Then a.a.s.\ there is a constant $I>0$, such that after $I$ rounds of the parallel peeling algorithm are applied to $H$, every component of the remaining graph, denoted by $H_I$, has size $O(\log n)$.
\end{prop}

\no {\bf Proof of Theorem~\ref{thm:main}.\ } Let $I$ be a constant chosen to satisfy Proposition~\ref{prop} and let $H_I$ be the remaining graph after $I$ rounds of the parallel peeling algorithm. Then, a.a.s.\ every component of $H_I$ contains $O(\log n)$ vertices. By Lemma~\ref{lem:small}, we may assume that each component has average degree at most $r/(r-1)+\eps$ for any constant $\eps>0$. Take an arbitrary constant $C$ of $H_I$. Let $C_0,C_1,\ldots,$ denote the process produced by running the parallel peeling algorithm on $C_0=C$. By Lemma~\ref{lem:small}, we may assume that each $C_i$ has average degree at most $r/(r-1)+\eps$. Let $\rho_i$ denote the proportion of vertices in $C_i$ with degree at least $k$. Then $k\rho_i\le r/(r-1)+\eps$ for every $i\ge 0$; i.e.\ $\rho_i\le \rho:=r/k(r-1)+\eps/k$. By our assumption on $k$ and $r$, we always have $\rho<1$. Since all vertices with degree less than $k$ are removed in each step of the algorithm, we have
$|V(C_{i+1})|\le \rho |V(C_i)|$ for every $i\ge 0$. This immediately gives $s(C)\le (\log\log n+O(1))/\log\rho^{-1}$. Since $\eps>0$ can be taken arbitrarily small, we have $s(C)\le (a^*_{r,k}+o(1))\log\log n$. This holds a.a.s.\ for every component of $H_I$. Hence, a.a.s.\ $s(\H_r(n,c/n^{r-1}))\le I+(a^*_{r,k}+o(1))\log\log n=(a^*_{r,k}+o(1))\log\log n$. \qed

\end{document}